\documentclass{amsart}

\usepackage{tikz}

\newtheorem{theorem}{Theorem}
\newcommand{\bt}{\begin{theorem}}
\newcommand{\et}{\end{theorem}}
\newtheorem{lemma}{Lemma}
\newcommand{\bl}{\begin{lemma}}
\newcommand{\el}{\end{lemma}}
\newtheorem{corollary}{Corollary}
\newcommand{\bc}{\begin{corollary}}
\newcommand{\ec}{\end{corollary}}
\newcommand{\beq}{\begin{equation}}
\newcommand{\eeq}{\end{equation}}
\newcommand{\benum}{\begin{enumerate}}
\newcommand{\eenum}{\end{enumerate}}

\newcommand{\R}{\ensuremath{\mathbf R}}

\newcommand{\mcw}{\ensuremath{ \mathcal W}}

\newcommand{\mba}{\ensuremath{ \mathbf a}}
\newcommand{\mbb}{\ensuremath{ \mathbf b}}

\newcommand{\mbe}{\ensuremath{ \mathbf e}}

\newcommand{\mbu}{\ensuremath{ \mathbf u}}
\newcommand{\mbv}{\ensuremath{ \mathbf v}}
\newcommand{\mbw}{\ensuremath{ \mathbf w}}
\newcommand{\mbx}{\ensuremath{ \mathbf x}}
\newcommand{\mby}{\ensuremath{ \mathbf y}}
\newcommand{\mbz}{\ensuremath{ \mathbf z}}

\newcommand{\Rm}{\ensuremath{ \mathbf{R}^m }}

\newcommand{\Rn}{\ensuremath{ \mathbf{R}^n }}

\newcommand{\bsmallmat}{\left(\begin{smallmatrix}}
\newcommand{\esmallmat}{\end{smallmatrix}\right)}

\DeclareMathOperator{\cone}{\text{cone}}
\DeclareMathOperator{\conv}{\text{conv}}

\newcommand{\bmat}{\left(\begin{matrix}}
\newcommand{\emat}{\end{matrix}\right)}

\DeclareMathOperator{\qqand}{\qquad\text{and}\qquad}

\DeclareMathOperator{\vectorvn}{\left( \begin{matrix} v_1 \\ \vdots \\ v_n \end{matrix}\right)}

\DeclareMathOperator{\vectorxn}{\left( \begin{matrix} x_1 \\ \vdots \\ x_n \end{matrix}\right)}

\DeclareMathOperator{\vectorsmallan}{\left( \begin{smallmatrix} a_1 \\ \vdots \\ a_n \end{smallmatrix}\right)}

\DeclareMathOperator{\vectorsmallun}{\left( \begin{smallmatrix} u_1 \\ \vdots \\ u_n \end{smallmatrix}\right)}
\DeclareMathOperator{\vectorsmallvn}{\left( \begin{smallmatrix} v_1 \\ \vdots \\ v_n \end{smallmatrix}\right)}

\DeclareMathOperator{\vectorsmallxn}{\left( \begin{smallmatrix} x_1 \\ \vdots \\ x_n \end{smallmatrix}\right)}

\DeclareMathOperator{\vectorsmallbm}{\left( \begin{smallmatrix} b_1 \\ \vdots \\ b_m \end{smallmatrix}\right)}

\DeclareMathOperator{\vectorsmallym}{\left( \begin{smallmatrix} y_1 \\ \vdots \\ y_m \end{smallmatrix}\right)}

\newcommand{\mbo}{\ensuremath{\mathbf 0}}
\DeclareMathOperator{\col}{\text{col}}

\newtheorem{problem}{Problem}
\newcommand{\bprob}{\begin{problem}}
\newcommand{\eprob}{\end{problem}}

\usepackage{amssymb,latexsym}

\title[Polyhedra and the Farkas lemma]{Polytopes, polyhedra, and the Farkas lemma}       

\author{Melvyn B. Nathanson}

\address{Department of Mathematics\\Lehman College (CUNY)\\Bronx, NY 10468}

\email{melvyn.nathanson@lehman.cuny.edu}



\thanks{Supported in part by  the PSC-CUNY Research Award Program, grant \#63117-00 51.}

\date{\today}

\begin{document}

\begin{abstract}
The Farkas lemma is proved and then applied to obtain a structure theorem for polyhedra.  
These notes are based on a talk in the New York Number Theory Seminar on October, 20, 2022.  
\end{abstract}
\maketitle

\section{Notation} 

For vectors $\mbu = \vectorsmallun$ and $\mbv = \vectorsmallvn$ in \Rn\ 
we write $\mbu \leq \mbv$ (and also  $\mbv \geq \mbu$) 
if $u_i \leq v_i$ for all $i \in \{1,2,\ldots, n\}$. 
If $\mbu \leq \mbv$ and $t \geq 0$, then $t\mbu \leq t\mbv$. 
We denote both the \emph{zero vector} in \Rn\ and the zero $m \times n$ matrix 
in $\R^{m,n}$ by $\mbo$. 
The vector $\mbv$ is  \emph{nonnegative}\index{nonnegative vector}  
if $\mbv \geq \mbo$ and  \emph{nonpositive}\index{nonpositive vector}  
if $\mbv \leq \mbo$. 
The transpose of the column vector $\mbu = \vectorsmallun$  is the row vector 
$\mbu^t = (u_1, \ldots, u_n)$ and   
\[
\mbu^t \mbv = (u_1, \ldots, u_n)  \vectorvn = \sum_{i=1}^n u_iv_i. 
\] 
If  $\mbu \leq \mbo$ and $\mbv \geq \mbo$, then $\mbu^t \mbv  \leq 0$. 

Let $A^t$ denote the transpose of the matrix $A$.
The matrix $A = \bmat a_{i,j}\emat$ is  \emph{nonnegative}\index{nonnegative matrix}, 
denoted $A \geq \mbo$, 
if $a_{i,j} \geq 0$ for all $i$ and $j$, and \emph{nonpositive}\index{nonpositive matrix}, 
denoted $A \leq \mbo$, 
if $a_{i,j} \leq 0$ for all $i$ and $j$. 
The inequalities $A \leq \mbo$ and $\mbv \geq \mbo$ imply $ A\mbv \leq \mbo$.

\section{Convex sets and polyhedra} 

 A subset $X$ of a real vector space $V$ is \emph{convex}\index{convex set} 
 if $\lambda_1\mbx_1+\lambda_2\mbx_2 \in X$ 
 for all vectors $\mbx_1, \mbx_2 \in X$ and all nonnegative 
 numbers $\lambda_1, \lambda_2$ such that $\lambda_1+\lambda_2=1$. 
 The intersection of convex subsets of $V$ is a convex set in $V$. 
  For example, for every positive integer $k$, the set   
\[
\Delta_k = \left\{ \mathbf{\lambda} = \bmat \lambda_1 \\ \vdots \\ \lambda_k \emat  \in \R^k: 
\mathbf{\lambda} \geq \mbo \text{ and }\sum_{i=1}^k \lambda_i = 1 \right\}
\] 
is a convex subset of $\R^k$.  
The set $X$ is convex if and only if $\sum_{i=1}^k \lambda_i \mbx_i \in X$ 
for all $k \geq 1$, $(\mbx_1,\ldots, \mbx_k) \in X^k$, 
and $\bsmallmat \lambda_1 \\ \vdots \\ \lambda_k \esmallmat \in \Delta_k$. 
We call $\sum_{i=1}^k \lambda_i \mbx_i \in X$ a \emph{convex combination}\index{convex combination} 
of the $k$-tuple $(\mbx_1,\ldots, \mbx_k)$. 

The intersection of convex subsets of a real vector space $V$ is convex. 
The \emph{convex hull}\index{convex hull} of a  nonempty subset $W$ 
of $V$ is the intersection of all convex subsets of $V$ that contain $W$.  
This is the set $X = \conv(W)$ of all convex combinations of $k$-tuples of elements of $W$ 
for all $k \geq 1$.  
For $j \in \{1,2,\ldots, k\}$, the standard basis vector  $\mbe_j \in \R^k$ is the vector 
whose $j$th coordinate is 1 and 
whose $i$th coordinate is 0 for all $i \in \{1,2,\ldots, k\} \setminus \{ j\}$. 
The set $\Delta_k$ is the convex hull of the set $\{ \mbe_1,\ldots, \mbe_k\}$.  
We define $\conv(\emptyset) = \{\mbo\}$.

The \emph{closed half-space}\index{closed half-space}\index{half-space!closed} 
in \Rn\ defined by the  vector $\mba = \vectorsmallan \in \Rn$ and the scalar $b \in \R$ 
is the set 
\[
H(\mba,b) = \left\{ \mbx = \vectorsmallxn \in \Rn: \mba^t \mbx = \sum_{j=1}^n a_j x_j \leq b \right\}. 
\]
If $\mbx, \mby  \in H(\mba,b)$ and if $\lambda_1+\lambda_2 = 1$, then 
\[
\mba^t \left( \lambda_1 \mbx + \lambda_2 \mby \right)  =
\lambda_1 \mba^t \mbx + \lambda_2 \mba^t  \mby  \leq \lambda_1 b + \lambda_2 b = b 
\] 
and so $ \lambda_1 \mbx + \lambda_2 \mby \in H(\mba,b)$.  
Thus, every closed half-space is convex.  

If $\mba \neq 0$, then $H(\mba,b) \neq \{\mbo\} $ and $H(\mba,b) \neq   \Rn$.
We have 
$H(\mbo,b) = \emptyset$ if $b < 0$ and $H(\mbo,b) = \Rn$ if $b \geq 0$.

A \emph{polyhedron}\index{polyhedron} in \Rn\ is the intersection 
of a finite number of closed half-spaces, that is, a set $P$ of vectors  
$\mbx = \vectorsmallxn$ whose coordinates satisfy 
a finite number $m$ of linear inequalities of the form $\sum_{j=1}^n a_{i,j} x_j \leq b_i$.   
These inequalities define the  $m \times n$ matrix $A = \bmat a_{i,j} \emat$ 
and the vector $\mbb = \vectorsmallbm \in \R^m$,  which  
generate the polyhedron $P$ as follows:  
\[
P = \left\{\mbx \in \R^n: A \mbx \leq \mbb  \right\} = \bigcap_{i=1}^m H(\mba_i, b_i)  
\] 
where $\mba_i^t$ is the $i$th row vector of $A$. 
If $\mbx_1, \mbx_2 \in P$ and if $\lambda_1 + \lambda_2 = 1$ 
with $\lambda_1,  \lambda_2 \in \R^2_{\geq 0}$, then 
\[
A\left( \lambda_1\mbx_1 + \lambda_2 \mbx_2\right) =  \lambda_1 A \mbx_1 +   \lambda_2 A \mbx_2  
\leq  \lambda_1 \mbb +   \lambda_2 \mbb = \mbb
\]
and so every polyhedron is convex.

The inequality $\sum_{j=1}^n a_j x_j \geq b$ is equivalent to the inequality 
$\sum_{j=1}^n (-a_j) x_j \leq -b$, and the equation $\sum_{j=1}^n a_j x_j = b$  
is equivalent to the two inequalities $\sum_{j=1}^n a_j x_j \leq b$ 
and $\sum_{j=1}^n (-a_j) x_j \leq -b$.  
Thus, the set of vectors $\mbx \in \Rn$ that satisfy a finite number of 
inequalities of the form 
$\sum_{j=1}^n a_j x_j \leq b$, a finite number of 
inequalities of the form 
$\sum_{j=1}^n a_j x_j \geq b$, and   a finite number of equations  
$\sum_{j=1}^n a_j x_j = b$ is a polyhedron.

\section{Projections}

For all $n \geq 2$ and $k \in \{1,\ldots, n-1\}$, the \index{projection}\emph{projection} 
of $\R^n$ onto $\R^{n-k}$  is the llinear transformation $\pi_{n,k}:\Rn \rightarrow \R^{n-k}$ defined by 
\[
\pi_{n,k} \bsmallmat x_1 \\ \vdots \\ x_k \\ x_{k+1} \\ \vdots \\ x_n \esmallmat = \bsmallmat x_{k+1} \\ \vdots \\ x_n \esmallmat. 
\]
We have 
\beq                    \label{Farkas:proj}
\pi_{n,k} = \pi_{n-k+1,1} \cdots \pi_{n-1, 1}\pi_{n, 1}.
\eeq
For example,
\[
\pi_{4,1}\pi_{5,1} \bmat x_1\\ x_2\\ x_3 \\ x_4 \\ x_5 \emat = \pi_{4,1} \bmat  x_2\\ x_3 \\ x_4 \\ x_5 \emat 
=  \bmat x_3 \\ x_4 \\ x_5 \emat =  \pi_{5,2} \bmat x_1\\ x_2\\ x_3 \\ x_4 \\ x_5 \emat
\]

\bl       \label{Farkas:lemma:projection-convex}
Let $V$ and $V'$ be real vector spaces and let $T:V \rightarrow V'$ be a linear transformation.  
If $X$ is a convex subset of $V$, then $Y = T(X)$  is a convex subset of $V'$.  
In particular, the projection of a convex set in \Rn\ is a convex set.
\el

\begin{proof}
Let $\mby_1, \mby_2 \in Y$ 
and $(\lambda_1, \lambda_2) \in \Delta_2$.    
There exist $\mbx_1,\mbx_2 \in X$ such that $T(\mbx_1) = \mby_1$ and  $T(\mbx_2) = \mby_2$.  
If $X$ is convex, then $\lambda_1 \mbx_1 + \lambda_2 \mbx_2  \in X$ and so 
\[
\lambda_1 \mby_1 + \lambda_2 \mby_2 = \lambda_1 T\left (\mbx_1 \right)  + \lambda_2 T\left (\mbx_2   \right) 
= T\left ( \lambda_1 \mbx_1 + \lambda_2 \mbx_2 \right) \in T(X) = Y.  
\]
Thus,  $Y$ is convex. 

Every projection $\pi_{n,k}:\Rn \rightarrow \R^{n-k}$ is a linear transformation 
and so  the projection of a convex set is convex. 
This completes the proof. 
\end{proof}

\bt        \label{Farkas:theorem:projection}
The projection of a polyhedron is a polyhedron.
\et

\begin{proof} 
Identity~\eqref{Farkas:proj} implies that it  suffices to prove the Theorem for $n \geq 2$ and $k=1$.  
The  projection  $\pi_{n,1}:\Rn \rightarrow \R^{n-1}$ is defined by 
\[
\pi_{n,1} \bsmallmat x_1 \\ x_2   \\ \vdots \\ x_n \esmallmat = \bsmallmat x_2  \\ \vdots \\ x_n \esmallmat.
\]

Let $P = \{\mbx \in \Rn: A\mbx \leq \mbb\}$ be the polyhedron in \Rn\ 
defined by the $m \times n$ matrix $A = \bmat a_{i,j} \emat$ 
and the vector $\mbb = \vectorsmallbm \in \Rm$. 
We have $\mbx = \vectorsmallxn \in P$ if and only if the following $m$ inequalities are satisfied:
\begin{align*}
a_{1,1}x_1 + a_{1,2} x_2 + \cdots  + a_{1,n}x_n & \leq b_1 \\
a_{2,1}x_1 + a_{2,2} x_2 + \cdots  + a_{2,n}x_n & \leq b_2 \\
\vdots & \\ 
a_{m,1}x_1 + a_{m,2} x_2 + \cdots  + a_{m,n}x_n & \leq b_m
\end{align*}
Partition the set $S = \{1, 2, \ldots, m\}$ as follows: 
\begin{align*} 
S_1 & = \{p \in S: a_{p,1} > 0 \} \\
S_{-1} & = \{q \in S: a_{q,1} < 0 \} \\ 
S_0 & = \{r \in S: a_{r,1} = 0 \}. 
\end{align*}
Let $\mbx = \vectorsmallxn \in P$.  If $p \in S_1$, then  
\[
x_1  \leq  -\frac{a_{p,2}}{a_{p,1}} x_2 - \cdots   -\frac{a_{p,n}}{a_{p,1}} x_n +  \frac{b_p}{a_{p,1}}. 
\]
If $q \in S_{-1}$, then  
\[
  -\frac{a_{q,2}}{a_{q,1}} x_2 - \cdots   -\frac{a_{q,n}}{a_{q,1}} x_n +  \frac{b_q}{a_{q,1}} \leq x_1.
\]
For all $p \in S_1$ and $q \in S_{-1}$ we have  
\[
 -\frac{a_{q,2}}{a_{q,1}} x_2 - \cdots   -\frac{a_{q,n}}{a_{q,1}} x_n +  \frac{b_q}{a_{q,1}} \leq x_1
  \leq  -\frac{a_{p,2}}{a_{p,1}} x_2 - \cdots   -\frac{a_{p,n}}{a_{p,1}} x_n +  \frac{b_p}{a_{p,1}} 
\]
and so 
\beq              \label{Farkas:proj-1}
\left( \frac{a_{p,2}}{a_{p,1}}  -\frac{a_{q,2}}{a_{q,1}} \right) x_2 + \cdots  
 + \left( \frac{a_{p,n}}{a_{p,1}}  - \frac{a_{q,n}}{a_{q,1}} \right) x_n  \leq  \frac{b_p}{a_{p,1}} - \frac{b_q}{a_{q,1}}. 
\eeq
For all $r \in S_0$ we have 
\beq              \label{Farkas:proj-0}
a_{r,2} x_2 + \cdots  + a_{r,n}x_n  \leq b_r.
\eeq
Thus, if $\mbx \in P$, then the vector $\pi_{n,1}(\mbx)$ is in the polyhedron 
defined by the linear  inequalities~\eqref{Farkas:proj-1} 
and~\eqref{Farkas:proj-0} for all $p,q,$ and $r$. 
There are $|S_1| \ |S_{-1}| + |S_0|$  such inequalities.

Conversely, if the vector $\bsmallmat x_2 \\ \vdots \\ x_n \esmallmat \in \R^{n-1}$ satisfies 
inequalities~\eqref{Farkas:proj-1} 
and~\eqref{Farkas:proj-0} for all $p,q,$ and $r$, then 
\beq              \label{Farkas:proj-2}
 -\frac{a_{q,2}}{a_{q,1}} x_2 - \cdots   -\frac{a_{q,n}}{a_{q,1}} x_n +  \frac{b_q}{a_{q,1}}
  \leq  -\frac{a_{p,2}}{a_{p,1}} x_2 - \cdots   -\frac{a_{p,n}}{a_{p,1}} x_n +  \frac{b_p}{a_{p,1}} 
\eeq 
for all $p$ and $q$, and so 
\begin{align*}
\max_{q\in S_{-1}} & \left( -\frac{a_{q,2}}{a_{q,1}} x_2 - \cdots   -\frac{a_{q,n}}{a_{q,1}} x_n +  \frac{b_q}{a_{q,1}}  \right)  
 \leq \min_{p\in S_1} \left( -\frac{a_{p,2}}{a_{p,1}} x_2 - \cdots   -\frac{a_{p,n}}{a_{p,1}} x_n +  \frac{b_p}{a_{p,1}}\right). 
\end{align*}
The left side of this inequality is $-\infty$ if $S_{-1} = \emptyset$ and the right 
side of this inequality is $\infty$ if $S_1 = \emptyset$. 
If $x_1$ is any real number such that 
\begin{align*}
\max_{q\in S_{-1}} & \left( -\frac{a_{q,2}}{a_{q,1}} x_2 - \cdots   -\frac{a_{q,n}}{a_{q,1}} x_n +  \frac{b_q}{a_{q,1}}  \right) \\ 
&  \leq x_1  \leq \min_{p\in S_1} \left( -\frac{a_{p,2}}{a_{p,1}} x_2 - \cdots   -\frac{a_{p,n}}{a_{p,1}} x_n +  \frac{b_p}{a_{p,1}}\right) 
\end{align*}
then the vector $\mbx = \vectorsmallxn \in \Rn$ is in the polyhedron $P$, 
and $\pi_{n,1}(\mbx) = \bsmallmat x_2 \\ \vdots \\ x_n \esmallmat$. 
Thus, the projection of the polyhedron $P$ in \Rn\ is the polyhedron in $\R^{n-1}$ 
defined by the inequalities~\eqref{Farkas:proj-1} 
and~\eqref{Farkas:proj-0}. 
This completes the proof. 
\end{proof}

The method of proof of Theorem~\ref{Farkas:theorem:projection}  
is sometimes called \index{Fourier-Motzkin elimination}\emph{Fourier-Motzkin elimination}.

\section{Convex cones} 
A \emph{cone}\index{cone} in $\R^n$ is a nonempty set $C$ 
such that $\lambda \mbw \in C$ for all 
$\mbw \in C$ and $\lambda \geq 0$.  
If $\mbw \in C$, then $\mbo = 0\mbw \in C$ and so every cone contains the zero vector.  
The set $\{ \mbo \}$ in \Rn\ is a cone.  The cone $C$ is bounded if and only if $C = \{\mbo\}$. 

A \emph{convex cone}\index{convex cone}\index{cone!convex} 
in $\R^n$ is a cone that is convex.   
The set $\{ \mbo \}$ in \Rn\ is a convex cone.  
If the cone $C$ is convex, then  for all $\mbw, \mbw' \in C$ we have 
\[
\mbw + \mbw' = 2\left(\frac{1}{2} \mbw + \frac{1}{2} \mbw' \right)  \in C.
\]
Conversely, if $C$ is a cone such that $\mbw + \mbw' \in C$ 
for all $\mbw, \mbw' \in C$, then for all $(\lambda_1, \lambda_2) \in \Delta_2$ 
we have $\lambda_1\mbw, \lambda_2\mbw' \in C$ and so $\lambda_1\mbw + \lambda_2\mbw' \in C$.
Thus, the cone $C$ is convex if and only if $\mbw + \mbw' \in C$ 
for all $\mbw, \mbw' \in C$. 

The union of cones is a cone, but the union of convex cones is not necessarily convex.  
For example, the union of two distinct one-dimensional subspaces in \Rn\  
(that is, two lines passing through the origin)  is a cone that is not convex. 

The vector space $\R^{n^2}$ can be viewed as the set of $n \times n$ matrices.  
An $n \times n$ matrix $A$ is \emph{positive semidefinite}\index{positive semidefinite} 
if $\mbx^tA\mbx \geq 0$ for all $\mbx \in \Rn$. The set of positive semidefinite matrices 
is a convex cone in $\R^{n^2}$.

A \emph{conic combination}\index{conic combination} of a finite sequence 
of vectors  $ (\mbw_1,\ldots, \mbw_n)$  in $\R^n$ is a vector of the form 
$\lambda_1\mbw_1+\cdots + \lambda_n \mbw_n$ 
for some $(\lambda_1,\ldots, \lambda_n) \in \R^n_{\geq 0}$.  
A nonempty set  $\mathcal{W}$ of vectors in \Rn\ 
\emph{conically generates}\index{conically generate}\index{generate!conically} the cone $C$ 
if $C$ is the set of all conic combinations 
of finite sequences of elements of $\mathcal{W}$.  
We write $C = \cone(\mathcal{W})$ if $C$ is the cone generated by $\mathcal{W}$.   
A  cone   is \emph{finitely generated}\index{cone!finitely generated} 
if  it is the set of all conic combinations of a finite set of vectors. 
The zero cone $C = \{\mbo\}$ is the cone generated by the empty set of vectors in \Rn.

\bt         \label{Farkas:theorem:cone-W}
The cone generated by a nonempty set of vectors is convex. 
\et 

\begin{proof}
Let $\mathcal{W}$ be a nonempty set of vectors in \Rn\ and let $C = \cone(\mathcal{W})$. 
The set $C$ is nonempty because $\mathcal{W}$ is nonempty.  
A conic combination of conic combinations of finite sequences of vectors in $\mathcal{W}$ 
is also a conic combination of a finite sequence of vectors in $\mathcal{W}$, and so $C$ is a cone.  
 In particular, if $\mbw \in C$ and $\mbw' \in C$, then $\mbw + \mbw' \in C$. 
 Thus, $C$ is convex.  This completes the proof. 
\end{proof}

\bt        \label{Farkas:theorem:cone}
Let $W$ be an $n \times m$ matrix . 
The set 
\begin{align*}
C(W) =  \left\{ W\mathbf{y}: \mathbf{y} \in \R^m \text{ and }\mathbf{y} \geq \mbo \right\}
\end{align*}
is a convex cone in \Rn.  
\et 

\begin{proof}
For $j \in \{1,\ldots, m\}$, let $\mbw_j \in \Rn$ be the the $j$th column vector of the matrix $W$ 
and let $\mcw = \{\mbw_1,\ldots, \mbw_m\}$.   
We have 
\begin{align*}
C(W) & = \left\{ W\mathbf{y}: \mathbf{y} \in \R^m \text{ and }\mathbf{y} \geq \mbo \right\} \\
& = \left\{y_1\mbw_1+\cdots + y_m \mbw_m: y_j \geq 0 \text{ for all } j = 1,\ldots, m \right\} 
\end{align*} 
and so $C(W)$ is the set of all conic combinations of the sequence of vectors 
$(\mbw_1,\ldots, \mbw_m)$, that is, $C(W) = \cone(\mathcal{W})$.
By Theorem~\ref{Farkas:theorem:cone-W}, the cone $C(W)$ is convex.  
This completes the proof. 
\end{proof}

A \emph{polyhedral cone}\index{cone!polyhedral}\index{polyhedral cone} in $\R^n$ is a 
polyhedron that is a   cone. 

\bl                \label{Farkas:lemma:ConicalPolyhedron}
For every positive integer $m$ and every  $m \times n$ matrix $A$, the polyhedron 
\[
P(A) =  \left\{\mbx \in \R^n: A \mbx \leq \mbo  \right\}
\]
is a polyhedral cone.   
Every polyhedral cone in \Rn\ is of this form for some matrix $A$.  
\el

\begin{proof}
If $A\mbx  \leq \mbo$, then $A(\lambda \mbx) = \lambda A \mbx \leq \mbo$
for all $\lambda \geq 0$ and so the polyhedron 
$\left\{\mbx \in \R^n: A \mbx \leq \mbo  \right\}$ is a cone. 

Let $A = \bmat a_{i,j} \emat$ and $\mbb = \vectorsmallbm$. 
Suppose that the polyhedron $P = \left\{\mbx \in \R^n: A \mbx \leq \mbb  \right\}$ 
is a cone.  
We have $\mbo \in P$ and so $\mbo  = A\mbo \leq \mbb$. 
Therefore, $A\mbx \leq \mbo$ implies $A\mbx \leq \mbb$ 
and  $P$ contains the set $ \left\{\mbx \in \R^n: A \mbx \leq \mbo  \right\}$.

If $P$ contains a vector \mbx\ such that $A\mbx \not\leq \mbo$, 
then the $i$th coordinate of the vector $A\mbx$ is positive for some $i \in \{1,\ldots,  m\}$. 
This coordinate is $\sum_{j=1}^n a_{i,j} x_j $ and satisfies the inequality 
\[
0 < \sum_{j=1}^n a_{i,j} x_j \leq b_i.
\]
Because $P$ is a cone,  for all $\lambda  > 0$ we have $\lambda  \mbx \in P$  
and so $ \lambda  A\mbx = A(\lambda  \mbx) \leq \mbb$. 
The $i$th coordinate of the vector $\lambda  A\mbx$ satisfies the inequality 
\[
0 <  \lambda  \sum_{j=1}^n a_{i,j} x_j \leq b_i 
\]
for all $\lambda  > 0$, which is absurd.  Therefore, 
$P = \left\{\mbx \in \R^n: A \mbx \leq \mbo  \right\}$. 
This completes the proof. 
\end{proof}

\bl       \label{Farkas:lemma:projection-convex}
Let $V$ and $V'$ be real vector spaces and let $T:V \rightarrow V'$ be a linear transformation.  
If $X$ is a cone in $V$, then $Y = T(X)$  is a cone in $V'$.  
In particular, the projection of a cone in \Rn\ is a cone.
\el

\begin{proof}
Let $C$ be a cone in $V$.  If $\mbw \in T(C)$, then $\mbw = T(\mbv)$ for some $\mbv \in V$.
For all $\lambda \geq 0$ we have $\lambda \mbv \in C$ and so 
\[
\lambda \mbw = \lambda T(\mbv) = T(\lambda \mbv) \in T(C) 
\]
and so $T(C)$ is a cone. 

Every projection $\pi_{n,k}:\Rn \rightarrow \R^{n-k}$ is a linear transformation 
and so  the projection of a cone is a cone. 
This completes the proof.  
\end{proof}

\emph{Notation.} For $\mby = \bsmallmat y_1 \\ \vdots \\ y_p \esmallmat \in \R^p$ 
and $\mbx \in \vectorsmallxn \in \Rn$, 
let $\bmat \mby\\ \mbx \emat = \bsmallmat y_1 \\ \vdots \\ y_p \\ x_1 \\ \vdots \\ x_n \esmallmat  
\in \R^{p+n}$. 
Let $I_n$ be the $n \times n$ identity matrix and $\mbo_{p,n}$ the $p \times n$ zero matrix.

\bt[Weyl]           \label{Farkas:theorem:Weyl} 
Every finitely generated convex cone is a polyhedral cone.
\et

\begin{proof}
Let $C$ be the convex cone in \Rn\ generated by the finite set of vectors 
$\{\mbw_1, \ldots, \mbw_p\}$.  
Let  $W$ be the $n \times p$ matrix whose column vectors are $ (\mbw_1,\ldots, \mbw_p)$. 
Then     
\[
C = \left\{ W\mathbf{y}: \mathbf{y} \in \R^p \text{ and }\mathbf{y} \geq \mbo \right\}. 
\]
Consider the set 
\[
P = \left\{ \bmat \mby \\ W\mby \emat: \mby \in \R^p,  \mby \geq \mbo   \right\} 
\subseteq \R^{p+n}.
\]
We have
\begin{align*}
P & =  \left\{ \bmat \mby \\ \mbx \emat: \mby \in \R^p, \mby \geq \mbo, \mbx = W\mby  \right\}  \\
& =  \left\{ \bmat \mby \\ \mbx \emat: \mby \in \R^p, \mbx \in \R^n, -\mby \leq \mbo, 
W\mby - \mbx  \leq \mbo,   - W\mby + \mbx \leq \mbo  \right\} \\
& =  \left\{  \mbz \in \R^{p+n}: \tilde{W}  \mbz\leq \mbo \right\}
\end{align*} 
where $\tilde{W}$ is the $(p+2n)\times (p+n)$ matrix written in block form as 
\[
\tilde{W} = \bmat -I_p & \mbo_{p,n} \\ W & -I_n \\ -W & I_n \emat. 
\]  
It follows from Lemma~\ref{Farkas:lemma:ConicalPolyhedron} that $P$ is a polyhedral cone.  

Let $\mby \in \R^p$ and $\mbx \in \R^n$.  
The projection $\pi_{p+n,p}:\R^{p+n} \rightarrow \Rn$ is 
defined by $\pi_{p+n,p} \bmat \mby \\ \mbx \emat = \mbx$.
We have the cone 
\[
\pi_{p+n,p}(P) = \left\{  W\mby :  \mathbf{y} \in \R^p \text{ and } \mby \geq \mbo  \right\} = C.  
\]
By Theorem~\ref{Farkas:theorem:projection}, the set  
$C = \pi_{p+n,p}(P)$ is a polyhedron, and so $C$ is a polyhedral cone. 
This completes the proof. 
\end{proof}

The following result is an example of what is called a ``theorem of the alternative.''

\bt[Farkas lemma]            \label{Farkas:theorem:Farkas}
Let $W$ be an $n \times p$ matrix and let $\mbb \in \R^n$.  
Exactly one of the following two statements holds:
\benum
\item[(a)]
There  is a vector $\mby \in \R^p$ such that 
\[
\mby \geq \mbo \qqand W\mby = \mbb.
\] 
\item[(b)]
There  is a vector $\mbv \in \R^n$ such that 
\[
\mbv^tW \leq \mbo^t \qqand \mbv^t \mbb > 0.
\]
\eenum
\et

\begin{proof}
Suppose that alternative (a) holds and that $\mby \in \R^p$ satisfies 
$\mby \geq \mbo$ and $W\mby = \mbb$.  If $\mbv \in \R^n$ and 
$\mbv^t W \leq \mbo^t$, then 
\[
\mbv^t \mbb = \mbv^t \left( W\mby  \right)  = \left(\mbv^t  W \right) \mby  \leq 0 
\]
and so alternative (b) does not hold.

Let $\mbw_1, \ldots, \mbw_p$ be the column vectors of the matrix $W$   
and let   
\[
C(W) = \{W\mby: \mby \in \R^p \text{ and } \mby \geq \mbo\}
\]
be the convex cone 
in \Rn\ generated by the finite sequence $(\mbw_1, \ldots, \mbw_p)$.  
By Weyl's theorem (Theorem~\ref{Farkas:theorem:Weyl}), 
the finitely generated convex cone $C(W)$ is a polyhedral cone and so 
there is an $m \times n$ matrix $A$ such that 
\[
C(W) = P(A) = \{\mbx \in \Rn: A\mbx \leq \mbo\}. 
\]
If  alternative (a) does not hold, then $\mbb \notin C(W)$ 
and so $A\mbb \not\leq \mbo$. 
This means that the $i$th  coordinate of the vector $A\mbb$ is positive 
for some $i \in \{1,\ldots, n\}$.   Let $(a_{i,1},\ldots, a_{i,n})$  be the $i$th row of the matrix $A$ 
and let $\mbv = (a_{i,1},\ldots, a_{i,n})^t \in \Rn$. 
The $i$th coordinate of $A\mbb$ is $\mbv^t\mbb = \sum_{j=1}^n a_{i,j}b_j$ 
and so $\mbv^t\mbb > 0$. 

For all $j \in \{1,\ldots, p\}$, the $j$th column vector $\mbw_j = W\mbe_j $  
belongs to the convex cone $C(W)$, and so $A\mbw_j \leq \mbo$. 
The vector $A\mbw_j $ is the $j$th column vector of the matrix $AW$.  
Therefore, the $m \times h$ matrix $AW$ is nonnegative.  
Because $\mbv^t$ is the $i$th row of $A$, 
it follows that $\mbv^t W$ is the $i$th row of $AW$ and so  $\mbv^t W \leq \mbo^t$. 
Thus, alternative (b) holds if alternative (a) does not hold.  
This completes the proof of the Farkas lemma. 
\end{proof}

\section{A duality}
Associated to every $m \times n$ matrix $A$ are two convex sets in \Rn: 
The polyhedral cone 
\[
P = P(A) = \left\{ \mbx \in \Rn: A \mbx \leq \mbo \right\}
\]
and the finitely generated convex cone 
\[
C = C(A^t) = \left\{ A^t \mby: \mby \in \Rm \text{ and }  \mby \geq \mbo \right\}. 
\]
Consider, for example, the matrices 
\[
A_1 = \bmat -1 & 0 \\ 0 & -1 \\ -1 & 1 \emat \qqand 
  A_1^t = \bmat -1 & 0 & -1 \\ 0 & -1 & 1 \emat.
\]
We have the polyhedral cone 
\begin{align*}
P(A_1)  & = \left\{ \mbx \in \R^2: A \mbx \leq \mbo \right\} \\
& = \left\{ \bmat x_1 \\ x_2 \emat \in \R^2:  \bmat -1 & 0 \\ 0 & -1 \\ -1 & 1 \emat\bmat x_1 \\ x_2 \emat  
 \leq \bmat 0 \\ 0 \\ 0 \emat  \right\} \\
 & = \left\{  \bmat x_1 \\ x_2 \emat \in \R^2: \bmat -x_1 \\ -x_2 \\ x_2-x_1 \emat \leq \bmat 0 \\ 0 \\ 0 \emat  \right\} \\
 & = \left\{ \bmat x_1 \\ x_2 \emat \in \R^2: 0 \leq x_2 \leq x_1 \right\}
\end{align*}
and the convex cone 
\begin{align*}
C(A_1^t) & = \left\{ A^t \mby: \mby \in \R^3 \text{ and }   \mby \geq \mbo \right\} \\
& = \left\{ \bmat -1 & 0 & -1 \\ 0 & -1 & 1 \emat  \bmat y_1\\ y_2 \\ y_3 \emat  : y_1, y_2, y_3  \geq \mbo \right\} \\ 
& = \left\{\bmat -y_1 -y_3 \\ -y_2 + y_3  \emat : y_1, y_2, y_3  \geq \mbo \right\}. 
\end{align*}
Setting $x_1 = -y_1-y_3 \leq 0$ and $x_2 = -y_2+y_3 = -x_1 - y_1 - y_2 \leq -x_1$, 
we obtain  
\[
\left\{\bmat -y_1 -y_3 \\ -y_2 + y_3  \emat : y_1, y_2, y_3  \geq \mbo \right\} \\
\subseteq\left\{ \bmat x_1 \\ x_2 \emat \in \R^2: x_1 \leq 0 \text{ and }   x_2 \leq -x_1  \right\}.
\]
If $x_1 \leq 0$ and $x_2 \leq -x_1$, then choosing $y_1 = 0$, $y_2 = -x_1-x_2 \geq 0$, 
and $y_3 = -x_1 \geq 0$ gives 
\[
\bmat x_1 \\ x_2  \emat = \bmat -y_1 -y_3 \\ -y_2 + y_3  \emat  \in C(A_1^t)
\]
and so 
\[
C(A_1^t) =  \left\{ \bmat x_1 \\ x_2 \emat \in \R^2: x_1 \leq 0 \text{ and }  x_2 \leq -x_1 \right\}.
\]

  \begin{figure}     
 \begin{tikzpicture}[scale  = 0.8]
 \draw   [fill = red]  (0,0)--(8,0) -- (8,4)--(4,4)--(0,0);
  \draw   [fill = blue]  (0,0)--(-4,4) -- (-4,-4)--(0,-4)--(0,0);
  \end{tikzpicture} 
\caption{The cone $C(A_1^t)$ (blue) and the polyhedron $P(A_1)$ (red)}
\end{figure}
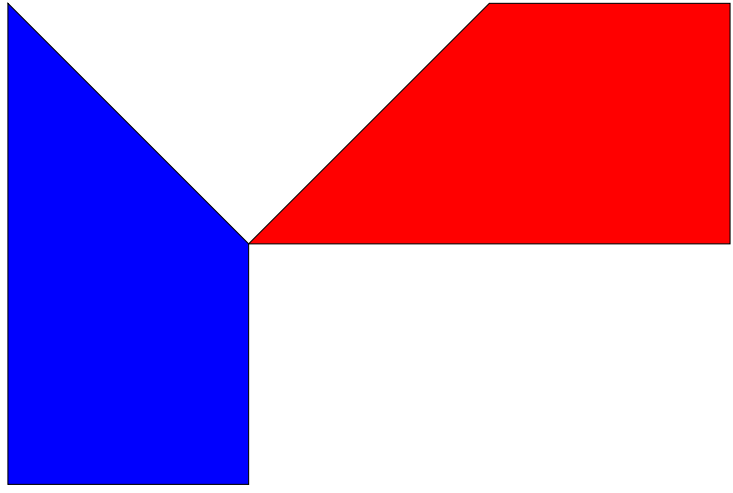

By Weyl's theorem, the cone $C(A_1^t)$ is a polyhedral cone.  
Indeed, $C(A_1^t)$ is defined by the inequalities $x_1 \leq 0$ and $x_2 \leq -x_1$. 
Consider the $2\times 2$ matrices  
\[
B_1 = \bmat 1 & 0 \\ 1 & 1 \emat \qqand B_1^t = \bmat 1 & 1 \\ 0 & 1 \emat.
\]
We have the polyhedral cone 
\begin{align*} 
P\left( B_1 \right) & = \left\{ \mbx \in \R^2: B_1\mbx \leq \mbo\right\} \\
& = \left\{  \bmat x_1\\ x_2 \emat \in \R^2: \bmat 1 & 0 \\ 1 & 1 \emat \bmat x_1\\ x_2 \emat \leq \bmat 0 \\ 0 \emat \right\} \\
& = \left\{  \bmat x_1\\ x_2 \emat \in \R^2: x_1 \leq 0 \text{ and }  x_1 + x_2 \leq 0 \right\} \\
& = C(A_i^t).
\end{align*} 
Thus,
\beq                 \label{Farkas:WAx}
\left\{ \mbx \in \R^2: B_1 \mbx \leq \mbo\right\} =  \left\{ A_1^t \mby: \mby \in \R^3 \text{ and }   \mby \geq \mbo \right\}. 
\eeq
The $2\times 2$ matrix  $B_1^t$ defines the convex cone 
\begin{align*}
 C(B_1^t) & = \left\{B_1^t \mby: \mby \in\R^2 \text{and } \mby \geq \mbo  \right\} \\
 & = \left\{ \bmat 1 & 1 \\ 0 & 1 \emat  \bmat y_1 \\ y_2 \emat: y_1,y_2 \geq 0  \right\} \\
  & = \left\{ \bmat y _1 + y_2 \\ y_2 \emat  \bmat y_1 \\ y_2 \emat: y_1,y_2 \geq 0  \right\} \\
  & = \left\{ \bmat x_1 \\ x_2 \emat \in \R^2: 0 \leq x_2 \leq x_1  \right\} \\
& = P(A_1)
\end{align*}
and so 
\beq                 \label{Farkas:AWx}
\left\{ \mbx \in \R^2: A_1 \mbx \leq \mbo \right\} =  \left\{B_1^t \mby: \mby \in\R^2 \text{and } \mby \geq \mbo  \right\}.
\eeq

Relations~\eqref{Farkas:WAx} and~\eqref{Farkas:AWx} illustrate the following result. 

\bt         \label{Farkas:theorem:AW} 
Let $A$ be an $m \times n$ matrix and let $B$ be an $n\times p$ matrix.  If  
\beq                            \label{Farkas:AW}
P(A) = \{ \mbx \in \Rn :  A\mbx \leq \mbo\} 
= \{ B \mby: \mby \in \R^p \text{ and } \mby \geq \mbo \}  = C(B) 
\eeq
then 
\beq                            \label{Farkas:WA}
P(B^t) = \{ \mbx \in \Rn :  B^t \mbx \leq \mbo\}
 = \{ A^t \mby: \mby \in \R^m \text{ and } \mby \geq \mbo \} = C(A^t). 
\eeq
\et

\begin{proof}
For $j \in \{1,\ldots, h\}$,  the $j$th column vector of $B$ is 
\[
 \mbb_j = B\mbe_j \in  \{ B \mby: \mby \in \R^p \text{ and } \mby \geq \mbo \}
\]    
and so 
\[
\mbb_j \in \{ \mbx \in \Rn :  A\mbx \leq \mbo\} 
\]
by  equation~\eqref{Farkas:AW}.  
The vector inequalities $A\mbb_j \leq \mbo$ for all $j \in \{1,\ldots, p\}$ 
imply the matrix inequality $AB \leq \mbo$.   The transpose of the matrix $AB$ 
is also nonpositive: 
\[
B^tA^t = (AB)^t \leq \mbo. 
\]
It follows that for all $\mby \geq \mbo$ we have 
\[
B^t\left(A^t \mby \right) = \left(B^t A^t  \right) \mby  \leq \mbo 
\] 
and so 
\[
C(A^t) =  \{ A^t \mby: \mby \in \R^m \text{ and } \mby \geq \mbo \} 
 \subseteq \{ \mbx \in \Rn :  B^t \mbx \leq \mbo\} = P(B^t). 
\]

We must prove the reverse inclusion: $P(B^t) \subseteq C(A^t)$. 
Let  $\mbx \in \Rn$ with $B^t \mbx \leq \mbo$. 
If   $A^t\mby \neq \mbx$ for all $\mby \in \Rm$ with $\mby \geq \mbo$, 
then the Farkas Lemma (Theorem~\ref{Farkas:theorem:Farkas}) implies that there exists a vector 
$\mbv \in \Rn$ such that 
\beq                            \label{Farkas:Btx}
(A\mbv)^t = \mbv^t A^t \leq \mbo^t \qqand \mbv^t \mbx > 0. 
\eeq
The inequalities $(A\mbv)^t  \leq \mbo^t$ and  $A\mbv \leq \mbo$ are equivalent. 
Equation~\eqref{Farkas:AW} implies $\mbv = B\mby$ for some $\mby \geq \mbo$. 
The inequalities $\mby  \geq \mbo$ and  $\mby^t \geq \mbo^t$ are equivalent.   
From $\mby^t \geq \mbo^t$ and $B^t \mbx \leq \mbo$ we obtain 
\[
\mbv^t \mbx = (B\mby)^t \mbx = \left(\mby^t B^t\right)  \mbx  = \mby^t  \left( B^t  \mbx \right) \leq 0
\]
which contradicts~\eqref{Farkas:Btx}.   Therefore, 
\[
P(B^t) = \{ \mbx \in \Rn :  B^t \mbx \leq \mbo\} \subseteq 
\{ A^t \mby: \mby \in \R^m \text{ and } \mby \geq \mbo \} =C(A^t). 
\]
This completes the proof. 
\end{proof}

We can now prove the converse of Weyl's Theorem.

\bt[Minkowski]            \label{Farkas:theorem:Minkowski}
Every polyhedral cone is a finitely generated convex cone. 
\et

\begin{proof}
Let  $A$ be an $m \times n$ matrix.  
Consider the polyhedral cone $P(A) = \{ \mbx \in \Rn: A\mbx \leq \mbo \}$ 
and the finitely generated convex cone 
$C(A^t) = \{ A^t \mby: \mby \in \Rm \text{ and } \mby \geq \mbo\}$. 
By Weyl's theorem, the cone $C(A^t)$ is polyhedral and so 
there is an $n \times p$ matrix $B$ such that 
\[
C(A^t)= \{ A^t \mby: \mby \in \Rm \text{ and } \mby \geq \mbo\} 
= \{ \mbx \in \Rn: B^t \mbx \leq \mbo\}= P(B^t). 
\] 
 Theorem~\ref{Farkas:theorem:AW} implies that 
\[
P(A)   = \{\mbx \in \Rn: A\mbx \leq \mbo\} 
= \{ B \mby: \mby \in \R^p \text{ and } \mby \geq \mbo\} = C(B)
\]
and so $P(A)$ is a finitely generated convex cone. 
This completes the proof. 
\end{proof}

\section{Structure of polyhedra} 

A \emph{polytope}\index{polytope} in $\R^k$  is the convex hull of a finite set of points in $\R^k$. 
For example, the set $\Delta_k = \conv(\mbe_1,\ldots, \mbe_k)$ is a polytope. 
A triangle in $\R^k$ is a polytope that is the convex hull of three non-collinear points.  
The convex hull of the set 
$\left\{ \bsmallmat 0 \\ 0 \esmallmat, \bsmallmat 2 \\ 0 \esmallmat, \bsmallmat 0 \\ 2 \esmallmat \right\}$ 
is the triangle 
\[
T = \left\{ \bmat x_1 \\ x_2 \emat: 0 \leq x_1 \leq 2 \text{ and } 0 \leq x_2 \leq 2-x_1\right\}.
\]
The triangle $T$ is also convexly generated by the set 
$\left\{ \bsmallmat 0 \\ 0 \esmallmat, \bsmallmat 2 \\ 0 \esmallmat, 
 \bsmallmat 1 \\ 1 \esmallmat,  \bsmallmat 1 \\ 1/2 \esmallmat,\bsmallmat 0 \\ 2 \esmallmat \right\}$ 
 and by the set $\{ (x,0): 0 \leq x \leq 2 \} \cup \{(0,2)\}$.  

The \emph{sum}\index{sum} (sometimes called the 
\emph{Minkowski sum}\index{Minkowski sum}\index{sum!Minkowski}) 
of subsets $C$ and $Q$ in $\R^k$ is the set 
\[
C+Q = \{ \mbu + \mbv : \mbu \in C \text{ and } \mbv \in Q \}.
\] 
We shall prove that every polyhedron is the sum of a finitely generated convex cone and a polytope, 
and that, conversely, every sum of a finitely generated convex cone and a polytope is a polyhedron.  

Here are two examples.  
Let 
\[
A_1 = \bmat -1 & 0 \\ 0 & -1 \\ -1 & -1 \emat \in \R^{3,2} \qqand \mbb_1 = \bmat  0 \\ 0 \\ -1 \emat \in \R^3.  
\]
The polyhedron 
\begin{align*}
P_1 &
= \left\{ \mbx \in \R^2: A_1 \mbx \leq  \mbb_1 \right\} \\
& = \left\{ \bmat x_1 \\ x_2  \emat: x_1 \geq 0, \ x_2 \geq 0, \ x_1+x_2 \geq 1 \right\}. 
\end{align*}
 is the sum of the convex cone  
\begin{align*}
C_1 & = \left\{ I_2 \mbx: \mbx \geq \mbo \right\} 
= \left\{ \bmat x_1 \\ x_2  \emat: x_1 \geq 0, \ x_2 \geq 0 \right\}. 
\end{align*}
 and the polytope $Q_1 = \Delta_2$.

  \begin{figure}                 
 \begin{tikzpicture}
 \draw  [fill = red]  (0,1)-- (0,4) -- (6,4) -- (6,0) -- (1,0) -- (0,1);
\draw [ultra thick] (0,1) -- (0,4); 
\draw [ultra thick] (1,0) -- (6,0); 
\draw [ultra thick] (0,1) -- (1,0);
 \end{tikzpicture} 
\caption{The polyhedron $P_1 = Q_1 + C_1$}
\end{figure}
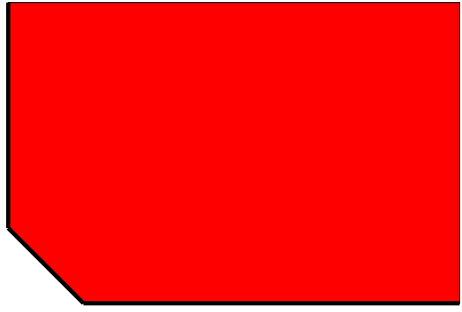

Let 
\[
A_2 = \bmat -1 & 1 \\ 1 & -1 \\ -1 & -1 \\ -2 & -1 \\ -1 & -2 \emat \in \R^{5,2} \qqand  
\mbb_1 = \bmat  4 \\ 4 \\ -3 \\ -4 \\ -4     \emat   \in \R^5.
\]
The polyhedron 
\begin{align*}
P_2 &
= \left\{ \mbx \in \R^2: A_2 \mbx \leq  \mbb_2 \right\} \\
& = \left\{ \bmat x_1 \\ x_2  \emat: -4 \leq  x_1 - x_2 \leq 4, \ x_1+x_2 \geq 3,  \ 
2x_1+x_2 \geq 4,  \ x_1+2x_2 \geq 4  \right\}. 
\end{align*}
is the sum of the convex cone  
\begin{align*}
C_2 & = \left\{ \bmat 1 \\ 1 \emat \bmat x_1  \emat: 
 x_1 \geq 0 \right\} = \left\{ \bmat x_1 \\ x_1  \emat: x_1 \geq 0 \right\}. 
\end{align*}
 and the polytope 
\begin{align*}
Q_2 & = \left\{ \lambda_1 \bmat 0 \\ 4 \emat +  \lambda_2 \bmat 1 \\ 2 \emat 
+  \lambda_3 \bmat 2 \\ 1 \emat +  \lambda_4 \bmat 4 \\ 0 \emat: 
\bmat \lambda_1 \\ \lambda_2 \\ \lambda_3 \\ \lambda_4 \emat \in \Delta_4  \right\}. 
\end{align*}

  \begin{figure}     
 \begin{tikzpicture}[scale  = 0.8]
 \draw   [fill = red]  (0,4)--(2,6) -- (6,6)-- (6,2) -- (4,0) -- (2,1) -- (1,2) -- (0,4);
 \draw [dashed] (0,4)-- (4,0);
\draw [ultra thick] (0,4)--(2,6); 
\draw [ultra thick] (1,2) -- (0,4); 
\draw [ultra thick]  (2,1) -- (1,2) ; 
\draw [ultra thick]   (2,1) -- (4,0); 
\draw [ultra thick] (4,0) -- (6,2);
 \draw  [fill = red]  (8,4) -- (12,0) -- (10,1) -- (9,2) -- (8,4);
 \end{tikzpicture} 
\caption{The polyhedron $P_2=Q_2+C_2$ and the polytope $Q_2$}
\end{figure}
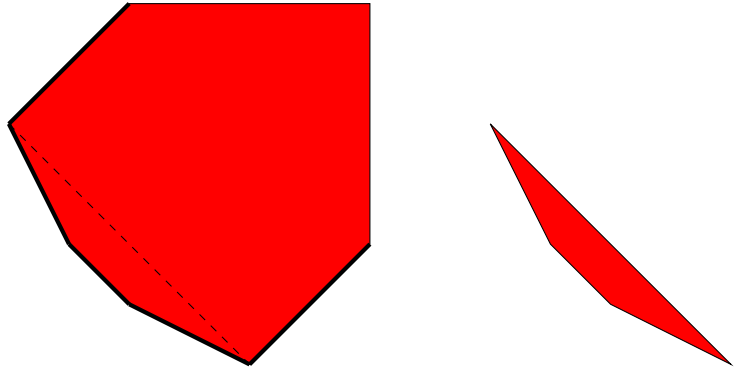

\bt             \label{Farkas:theorem:PQC}
A subset $P$ of \Rn\ is a polyhedron if and only if there is  a polytope $Q$ 
and a finitely generated convex cone $C$ such that $P = Q + C$.
\et

\begin{proof} 
Every vector in $\R^{n+1}$ can be written uniquely  in the form 
$\bmat \mbx \\ t \emat$ for some $\mbx \in \Rn$ and $t \in \R$. 
For every $\ell \times n$ matrix $A = \bmat a_{i,j} \emat$    
and vector $\mbb = \bsmallmat b_1 \\ \vdots \\ b_{\ell} \esmallmat \in \R^{\ell}$, 
let $\left( A|\mbb \right)$ be the ${\ell} \times (n+1)$ matrix whose $i$th row is 
$(a_{i,1}, \ldots , a_{i,n}, b_i )$ for all $i \in \{1,\ldots, {\ell}\}$. 
Every ${\ell} \times (n+1)$ matrix $\widehat{A}$ can be written uniquely in the form 
$\widehat{A} = \left( A|\mbb \right)$ for some ${\ell} \times n$ matrix $A$ and vector $\mbb \in \R^{\ell}$. 
For $\mbx = \vectorsmallxn \in \R^n$  and $t \in \R$, 
the $i$th coordinate of the vector $(A|\mbb) \bmat \mbx \\ t \emat \in \R^{\ell}$ is 
$\sum_{j=1}^n a_{i,j} x_j + b_i t$.   
Thus, the vector inequality $(A|\mbb) \bmat \mbx \\ t \emat \leq \mbo$ 
is equivalent to the ${\ell}$  inequalities  
\[
\sum_{j=1}^n a_{i,j} x_j \leq - tb_i 
\]
for all $i \in \{1,\ldots, {\ell} \}$ and so 
\beq                  \label{Farkas:Phat} 
(A|\mbb) \bmat \mbx \\ t \emat \leq \mbo \qquad\text{if and only if} \qquad A\mbx \leq -t\mbb.
\eeq

More generally, let $A_k = \bmat a_{i,j}^{(k)} \emat$ be an $m \times n_k$ 
matrix for all $k \in \{1,\ldots, \ell\}$ 
and let $n = \sum_{k=1}^{\ell} n_k$.  Let 
\[
A = (A_1| \cdots |A_k|\cdots |A_{\ell}) 
\]
be the $m \times n$ matrix obtained by concatenating the matrices $A_1, \ldots, A_{\ell}$.  
Thus, for $k \in \{1,\ldots, \ell\}$ and $j \in \{1,\ldots, n_k\}$,  
the $\left(\sum_{r=1}^{k-1} n_r + j \right)$th column of $A$ is 
\[
\col_j(A_k) = \bsmallmat a_{1,j}^{(k)}  \\ \vdots \\ a_{m,j}^{(k)} \esmallmat.
\]
The coordinate in row $i$ and column $\sum_{r=1}^{k-1} n_r + j$ 
of the matrix $A$ is  $a_{i,j}^{(k)}$.

Let $Q$ be a polytope in \Rn and let $C$ be a finitely generated convex cone in \Rn.    
We shall prove that the sumset $P = Q + C$ is a polyhedron.

Let $h \geq 2$ and let $\mbb_1,\ldots, \mbb_h$ be vectors in \Rn.  
For $k \in \{1,\ldots, h-1\}$, let $\widehat{P}$ be the convex cone in $\R^{n+1}$ 
conically generated by the finite sequence of vectors  
\[
\left( \bmat \mbb_1 \\ 1 \emat, \ldots, \bmat \mbb_k \\ 1 \emat, 
\bmat \mbb_{k+1} \\ 0 \emat, \ldots, \bmat \mbb_h \\ 0 \emat \right).
\]
Thus,
\begin{align*}
\widehat{P} 
& = \left\{ \sum_{j=1}^k \lambda_j \bmat \mbb_j \\ 1 \emat + \sum_{j=k+1}^h \lambda_j \bmat \mbb_j \\ 0 \emat   : \lambda_1, \ldots, \lambda_k, \lambda_{k+1}, \ldots, \lambda_h \geq 0 \right\}. 
\end{align*}

By Weyl's theorem (Theorem~\ref{Farkas:theorem:Weyl}) 
and by relation~\eqref{Farkas:Phat}, the  finitely generated convex 
cone $\widehat{P}$ is a polyhedral  cone, and so there is an $\ell \times n$ matrix $A$ and 
a vector $\mbb \in \R^{\ell}$  such that the  $\ell \times (n+1)$ 
matrix $(A|\mbb)$ satisfies  
\beq           \label{Farkas:Phat-2}
\widehat{P} = \left\{ \bmat \mbx \\ t \emat \in \R^{n+1}:  (A|\mbb) \bmat \mbx \\ t \emat \leq \mbo \right\} 
= \left\{ \bmat \mbx \\ t \emat \in \R^{n+1}: A\mbx \leq -t\mbb \right\}. 
\eeq

Let $(\mbw_1, \ldots, \mbw_k)$ be a finite  sequence of vectors  
that  convexly generates the polytope $Q$ 
and let  $(\mbw_{k+1}, \ldots, \mbw_h)$ be a  finite sequence of vectors 
that conically generates the cone $C$.   
If $\mbx \in Q + C$, then 
\[
\mbx = \mbx_Q + \mbx_C
\] 
with 
\[
\mbx_Q = \sum_{j=1}^k \lambda_j \mbw_j \in Q 
\qqand 
\mbx_C = \sum_{j=k+1}^h \lambda_j \mbw_j\in C
\]
and 
\[
\lambda_j \geq 0 \quad\text{for all $j \in \{1,\ldots,k,k+1,\ldots, h\}$} 
\qqand \sum_{j=1}^k \lambda_j = 1.  
\]
Thus,  
\begin{align*}
\bmat \mbx \\ 1\emat 
& = \bmat  \mbx_Q  \\ 1 \emat  + \bmat \mbx_C \\ 0\emat \\
& = \bmat  \sum_{j=1}^k \lambda_j \mbw_j  \\  \\  \sum_{j=1}^k \lambda_j\emat  
+ \bmat \sum_{j=k+1}^h \lambda_j \mbw_j \\  \\0  \emat  \\
& =   \sum_{j=1}^k \lambda_j  \bmat  \mbw_j   \\ 1 \emat  
+  \sum_{j=k+1}^h \lambda_j \bmat \mbw_j \\ 0 \emat  \\
& \in \widehat{P}.  
\end{align*} 
It follows from~\eqref{Farkas:Phat-2} that  $A\mbx \leq -\mbb$ 
and so $C+Q \subseteq \left\{ \mbx \in \Rn: A\mbx \leq -\mbb \right\}$. 

Similarly, if $\mbx \in \R^{n}$ and $A\mbx \leq -\mbb$, 
then  $\bmat \mbx \\ 1 \emat \in \widehat{P}$ and there exist 
nonnegative numbers $\lambda_1, \ldots, \lambda_k, \lambda_{k+1}, \ldots, \lambda_h$ 
such that 
\begin{align*}
\bmat \mbx \\ 1 \emat 
& = \sum_{j=1}^k \lambda_j \bmat \mbw_j \\ 1 \emat + \sum_{j=k+1}^h \lambda_j \bmat \mbw_j \\ 0 \emat \\
& =\bmat  \sum_{j=1}^k \lambda_j \mbw_j \\  \sum_{j=1}^k \lambda_j  \emat 
+\bmat  \sum_{j=k+1}^h \lambda_j  \mbw_j \\ 0 \emat. 
\end{align*} 
It follows that $ \sum_{j=1}^k \lambda_j  = 1$ and so 
\[
 \mbx_Q =  \sum_{j=1}^k \lambda_j \mbw_j \in Q \qqand  \mbx_C =  \sum_{j=k+1}^h \lambda_j  \mbw_j \in C.
\]
Thus, $\mbx = \mbx_Q + \mbx_C \in Q+C$ and  
$\left\{ \mbx \in \Rn: A\mbx \leq -\mbb \right\} \subseteq  Q+C $. 
This proves that 
\[
Q+C = \left\{ \mbx \in \Rn: A\mbx \leq -\mbb \right\}  
\] 
and so the sum of  a polytope and a finitely generated convex cone is a polyhedron.

Conversely, we shall  prove that every polyhedron is  the sum of a finitely generated 
convex cone and a polytope.  
Let $P$ be a polyhedron in \Rn.  There is an $m \times n$ matrix 
$A = \bmat a_{i,j}\emat$ and a vector $\mbb \in \Rm$ 
such that 
\[
P = \{ \mbx \in \Rn : A\mbx \leq \mbb\}. 
\]
Consider the $(m+1)\times (n+1)$ matrix 
\[
\widehat{A} = \bmat
a_{1,1} & a_{1,1} & \cdots & a_{1,n} & -b_1 \\
a_{2,1} & a_{2,2} & \cdots & a_{2,n} & -b_2 \\
\vdots & & & & \vdots \\
a_{m,1} & a_{m,2} & \cdots & a_{m,n} & -b_m \\
0 & 0 & \cdots &0 & -1
\emat
\]
By Lemma~\ref{Farkas:lemma:ConicalPolyhedron}, 
the set $\widehat{P}$ in $\R^{n+1}$ defined by 
\begin{align}     \label{Farkas:definePCQ}
\widehat{P} & = \left\{ \bmat \mbx \\ t \emat : \mbx \in \Rn, t \geq 0, A\mbx \leq t\mbb \right\} \\
& = \left\{ \bmat \mbx \\ t \emat \in \R^{n+1}: \widehat{A}  \bmat \mbx \\ t \emat  \leq \mbo \right\} 
\nonumber
\end{align}
is a polyhedral cone.  
We have $\mbx \in P$ if and only if $A\mbx \leq \mbb$ 
if and only if  $\bmat \mbx \\ 1 \emat \in \widehat{P}$ 
and so 
\[
P = \left\{ \mbx \in \Rn: \bmat \mbx \\ 1 \emat \in \widehat{P} \right\}. 
\]
By Minkowski's theorem (Theorem~\ref{Farkas:theorem:Minkowski}), the polyhedral 
cone $\widehat{P}$ is a finitely generated convex cone in $\R^{n+1}$, 
and so there is an $(n+1) \times h$ matrix $W$ such that 
\[
\widehat{P} = \{ W\mby: \mby \in \R^h \text{ and } \mby \geq 0\}. 
\]
The set $\widehat{P}$ is conically generated by the set of columns of $W$.  
For $j \in \{1,\ldots, h\}$, let $\bmat \mbw_j \\ t_j \emat$ be the $j$th column of $W$, 
where $\mbw_j \in \Rn$ and $t_j \in \R$.  
We have  $t_j \geq 0$ because  $\bmat \mbw_j \\ t_j \emat = W\mbe_j \in \widehat{P}$. 

Modify the matrix $W$ as follows:
If $t_j > 0$, then multiply the $j$th column of $W$ by $1/t_j > 0$.  If $t_j = 0$, do not 
change the $j$th column of $W$.  
We obtain a new matrix, which we shall also denote by $W$, whose $j$th column 
is of the form  $\bmat \mbw_j \\  \varepsilon_j \emat$ with $\varepsilon_j = 0$ or 1.
Because the columns of the original matrix $W$ have only been multiplied by positive numbers, 
the new matrix $W$ still conically generates $\widehat{P}$. 
Renumber the vectors $\mbw_1,\ldots, \mbw_h$ so that 
\[
\{\mbw_j: \varepsilon_j = 1 \} = \{\mbw_1,\ldots, \mbw_k\}
\qqand 
\{\mbw_j: \varepsilon_j = 0 \} = \{\mbw_{k+1},\ldots, \mbw_h\}. 
\]

Let  $Q$ be the polytope  in \Rn\ convexly generated by the set $\{\mbw_1,\ldots, \mbw_k\}$
and let $C$ be the convex cone in \Rn\ conically generated by the set $ \{\mbw_{k+1},\ldots, \mbw_h\}$. 
Let $Q = \{\mbo\}$ if $k=0$, that is, $\{\mbw_j: \varepsilon_j = 1 \} = \emptyset$,  
and let $C = \{\mbo\}$ if $k=h$. 

We have $\mbx \in P$ if and only if $\bmat \mbx \\ 1 \emat \in \widehat{P}$ if and only if 
there is a nonnegative vector $\mathbf{\lambda} = \bsmallmat \lambda_1 \\ \vdots \\ \lambda_{h} \esmallmat$ such that 
\begin{align*}
\bmat \mbx \\ 1 \emat 
& = W \mathbf{\lambda}    =  \sum_{j=1}^k \lambda_j \bmat \mbw_j \\ 1\emat +  \sum_{j=k+1}^h \lambda_j \bmat \mbw_j \\ 0 \emat   \\
& =  \bmat  \sum_{j=1}^k  \lambda_j \mbw_j \\  \\  \sum_{j=1}^k  \lambda_j \emat 
 +  \bmat \sum_{j=k+1}^h \lambda_j \mbw_j \\  \\ 0 \emat.
\end{align*}
It follows that $ \sum_{j=1}^k  \lambda_j=1$ and so 
 \[
 \sum_{j=1}^k  \lambda_j \mbw_j  \in Q 
\qqand 
 \sum_{j=k+1}^h \lambda_j \mbw_j  \in C. 
\]
We obtain  
\[
  \mbx = \sum_{j=1}^k  \lambda_j \mbw_j +  \sum_{j=k+1}^h \lambda_j \mbw_j \in Q+C. 
\]
Thus,  $P \subseteq Q+C$. 

Now we prove that $Q+C \subseteq P$.
If $\mbx \in Q+C$, then $\mbx = \mbx_Q + \mbx_C$, where  
  \[
\mbx_Q =  \sum_{j=1 }^k \lambda_j \mbw_j \in Q
 \qqand 
\mbx_C =   \sum_{ j=k+1  }^h \lambda_j \mbw_j \in C 
  \]
and $\mathbf{\lambda}  = \bsmallmat \lambda_1 \\ \vdots \\ \lambda_h \esmallmat \in \R^h$ 
 is a nonnegative vector such that $\sum_{j = 1}^k \lambda_j = 1$.  
Then 
\begin{align*}
\bmat \mbx \\ 1 \emat & =   \bmat \mbx_Q \\ 1 \emat + \bmat \mbx_C \\ 0 \emat \\
& =   \sum_{j=1 }^k \lambda_j  \bmat \mbw_j \\ 1 \emat + 
 \sum_{j=k+1}^h \lambda_j  \bmat \mbw_j \\ 0 \emat \\
 & \in \widehat{P}
  \end{align*}
  and so $\mbx \in P$.  Thus, $C+Q \subseteq P$.
It follows that $P = C+Q$ and so every polyhedron is the sum of a finitely generated 
convex cone and a polytope. 
 This completes the proof. 
\end{proof}

\section{Notes}
Much of the material in this paper derives from the classic paper of Weyl~\cite{weyl35}.  
The English translation is Weyl~\cite{weyl50}.  
Excellent introductions to convexity are the books by Eggleston~\cite{eggl58} 
and Lauritzen~\cite{laur13}.

\end{document}